\documentclass[prl,twocolumn, superscriptaddress]{revtex4}
\usepackage{graphicx}
\usepackage{amssymb}
\usepackage{amsmath}
\usepackage{ulem}
\usepackage{dcolumn}
\usepackage{epsfig}
\usepackage{bm}
\usepackage[latin1]{inputenc}
\usepackage{amsmath}
\usepackage[T1]{fontenc}
\usepackage[cm]{aeguill}
\usepackage{mathtools}
\usepackage{textcomp}
\usepackage{ulem}
\usepackage[english]{babel}
\usepackage{psfrag}
\usepackage{array}
\usepackage{eqnarray}
\usepackage{textcomp}
\usepackage{fancyhdr}
\usepackage{makeidx}
\def\Xint#1{\mathchoice
{\XXint\displaystyle \,  {#1}} 
{\XXint - \scriptstyle{#1}} 
{\XXint\scriptstyle\scriptscriptstyle{#1}} 
{\XXint\scriptscriptstyle\scriptscriptstyle{#1}} 
\!\int}
\def\XXint#1#2#3{{\setbox0=\hbox{$#1{#2#3}{\int}$ }
\vcenter{\hbox{$#2#3$ }}\kern-.54\wd0}}
\def\ddashint{\Xint=}
\def\dashint{\Xint-}
\newcommand{\ff}[1]{\mbox{\boldmath$#1$}}
\newcommand{\ffs}[1]{\mbox{\scriptsize \boldmath$#1$}}
\newcommand{\comment}[1]{}

\begin{document}
\title{Intermediate Asymptotics of the Capillary-Driven Thin Film Equation}

\author{M. Benzaquen}
\affiliation{Laboratoire de Physico-Chimie Th\'eorique, UMR CNRS Gulliver 7083, ESPCI ParisTech, PSL Research University, Paris, France}
\author{T. Salez}
\affiliation{Laboratoire de Physico-Chimie Th\'eorique, UMR CNRS Gulliver 7083, ESPCI ParisTech, PSL Research University, Paris, France}
\author{E. Rapha\"{e}l}
\affiliation{Laboratoire de Physico-Chimie Th\'eorique, UMR CNRS Gulliver 7083, ESPCI ParisTech, PSL Research University, Paris, France}

\date{\today}
\begin{abstract}
We present an analytical and numerical study of the two-dimensional capillary-driven thin film equation. In particular, we focus on the intermediate asymptotics of its solutions. Linearising the equation enables us to derive the associated Green's function and therefore obtain a complete set of solutions. Moreover, we show that the rescaled solution for any summable initial profile uniformly converges in time towards a universal self-similar attractor that is precisely the rescaled Green's function. Finally, a numerical study on compact-support initial profiles enables us to conjecture the extension of our results to the nonlinear equation.
\end{abstract}

\maketitle

Dimensional analysis is well understood through the Vaschy-Buckingham $\Pi$~theorem \cite{Buckingham1914}. Historically, this framework has already led to remarkably important results such as the expression of the hydrodynamical drag force on a sphere by Reynolds \cite{Reynolds1895}, the theory of turbulence by Kolmogorov \cite{Kolmogorov1942,Kolmogorov1941}, and the estimation of the nuclear explosion power by Taylor \cite{Taylor1950a,Taylor1950b}. Moreover, dimensional analysis is directly connected to the fundamental concepts of \textit{scaling} and \textit{self-similarity}, that appear in numerous situations such as fractals and diffusion \cite{Fourier1822}.
\smallskip

The powerful theory of \textit{intermediate asymptotics} developed in particular by Barenblatt goes one step further in understanding the deeper meaning of self-similarity \cite{Barenblatt1996}. In nonlinear problems, one may wonder what is the interest of finding exact particular solutions as there is no principle of superposition. Nonetheless, in certain cases, the self-similar solutions obtained for \textit{idealised} problems, or idealised initial conditions, represent the intermediate asymptotic regimes of the solutions of more general \textit{non-idealised} problems. Following Zeldovich in the foreword of \cite{Barenblatt1996}, one might even say that intermediate asymptotics is the key that somehow replaces the superposition principle in nonlinear physics. In other words, by paying the price of a loss of information at short times, one obtains a certain generality at intermediate times. Even so, such an asymptotic behaviour must be proved for any given initial condition, which often turns out to be a difficult task.
\smallskip

Furthermore, in the afterword of his book Barenblatt states \cite{Barenblatt1996}:
\textit{"However, there exist many problems of recognised importance where this technique has not yet been fully explained, but for which results of substantial value can be expected from its application."} One of these open problems is precisely the \textit{capillary-driven thin film equation} of interest \cite{Blossey2012,Craster2009,Oron1997}: 
\begin{equation}
\partial_T H+\partial_X\left(H^3\,\partial_X^{\,3} H\right)=0\ .
\label{TFE_AD}
\end{equation}
It governs the capillary evolution of the profile $H(X,T)$ of the free surface of a thin viscous liquid film: as soon as the curvature is nonconstant, this profile is unstable as the Laplace pressure drives a flow that is mediated by viscosity. Despite many efforts, this equation remains only partially solved. Nevertheless, in the past few years, several analytical \cite{Bowen2006,Myers1998,Salez2012} and numerical studies \cite{Bertozzi1998,Salez2012a} have been performed in order to gain insights into this mathematical problem. 
\begin{figure}[b!]
\includegraphics[width=1\columnwidth]{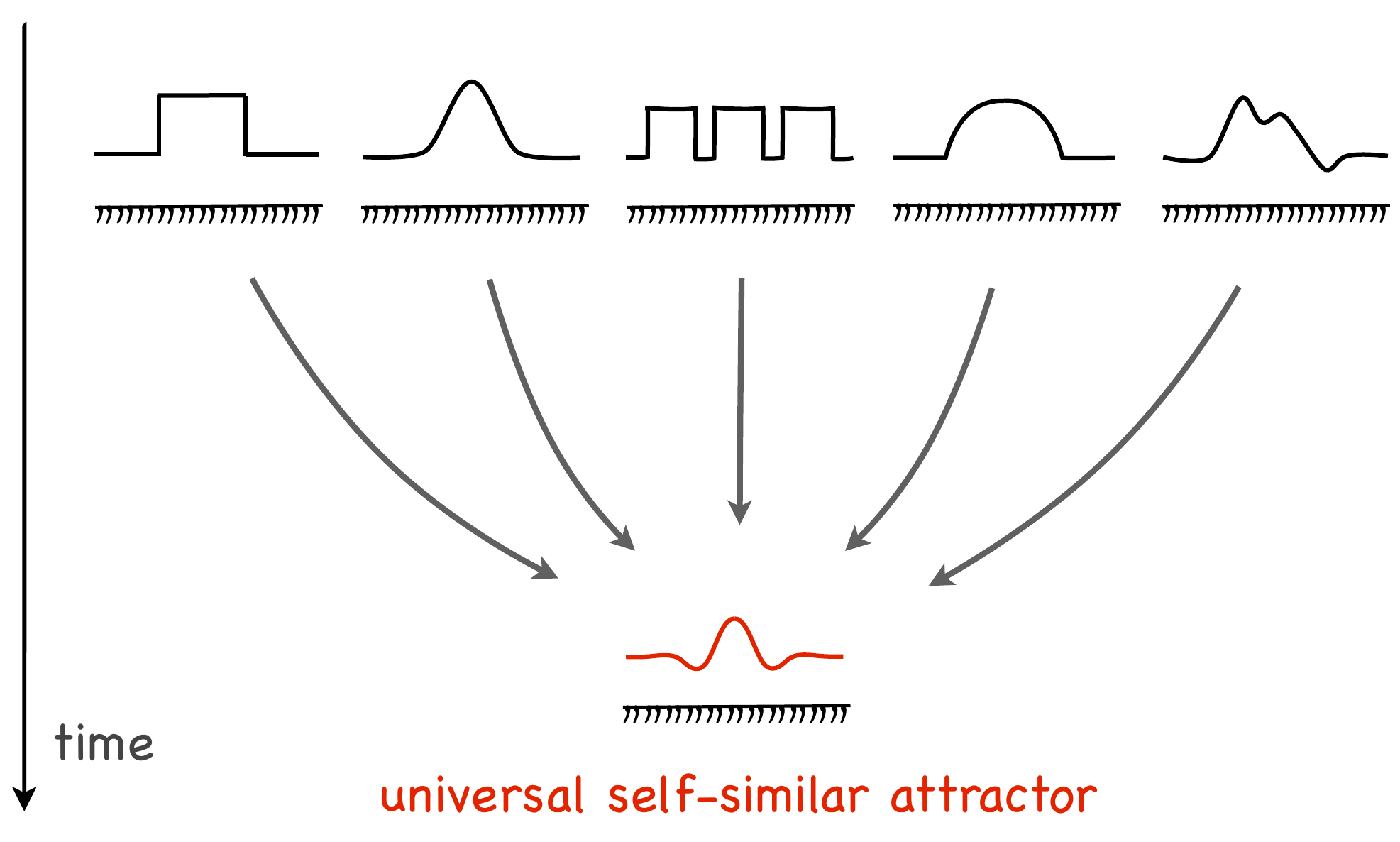}
\caption{Schematics of the intermediate asymptotics of the capillary-driven thin film equation (see Eq.~\eqref{TFE_AD}). No matter the initial condition, any summable profile converges in time towards a universal self-similar attractor.}
\label{Fig1}
\end{figure} 
\smallskip

It should be stressed that Eq.~\eqref{TFE_AD} is of tremendous importance in a variety of scientific fields such as polymer physics, physiology, biophysics, micro-electronics, surface chemistry, thermodynamics and hydrodynamics, since thin films are involved in modern mechanical and optical engineering processes, through lubrication, paints and coating. Gaining a complete understanding of these systems is a key step towards the development of molecular electronics, biomimetics, superadhesion and self-cleaning surfaces. Moreover, this equation may be crucial for understanding the nanorheology of ultra-thin polymer films, for which enhanced mobility effects have been predicted \cite{Brochard2000}, before being related to entanglements networks \cite{Si2005}, and observed in various experimental configurations such as: confinement \cite{Barbero2009,Bodiguel2006,Jones1999,Shin2007}, glassy state \cite {Fakhraai2008}, and dewetting onto slippery substrates \cite{Baumchen2009,Munch2005}. It may also govern the surface instabilities and pattern formation \cite{Amarandei2012,Closa2011,Mukherjee2011}, through the film preparation by spin-coating \cite{Baumchen2012,Raegen2010,Reiter2001,Stillwagon1990}. From all these examples, we understand the necessity of further exploring the solutions of the capillary-driven thin film equation. In particular, as for the diffusion equation \cite{Fourier1822}, it would be interesting to characterise the convergence of the solutions to some asymptotic self-similar regimes \cite{Stone2012}, as depicted in Fig.~\ref{Fig1}, as well as the existence of possible exotic self-similarities \cite{Sekimoto2012}. 
\smallskip

The present article is divided into three parts. In the first one, we recall the main ingredients of the physical model that describes two-dimensional capillary-driven flows. In the second one, we linearise the governing equation and we derive the Green's function that enables to calculate the general solution of the linear problem for any summable initial condition. In particular, we study the self-similar asymptotics of this solution (see Fig.~\ref{Fig1}). In the third part, we extend these ideas to the nonlinear equation through numerical solutions for compact-support initial profiles.
\newline

%%%%%%%%%%%%%%%%%%%%%%%%%%%%%%%%%%%%%%%%%%%%%%%%%%%%%%%%%%%%%%%%%%%%%%

\section{Physical model}
In this part, we describe the model and derive the capillary-driven thin film equation, within the lubrication approximation, for two-dimensional viscous flows. The evolution in time $t$ of a thin viscous film described by a profile of height $z=h(x,t)$ can be understood from the Laplace pressure $p(x,t)$, which arises due to curvature at the free interface. Considering small curvature gradients, one can write:
\begin{eqnarray}
p(x,t)&\simeq&-\gamma\,\partial_x^{\,2} h\ ,
\label{Laplace}
\end{eqnarray}
where  $\gamma$ is the liquid-air surface tension. In the lubrication approximation, the Stokes equation along the $x$ horizontal direction is given by \cite{Landau1987}:
\begin{eqnarray}
\partial_xp&=&\eta\,\partial_{z}^{\,2}v\ ,
\label{Stokes}
\end{eqnarray}
where $v(x,z,t)$ is the horizontal velocity and $\eta$ is the shear viscosity.
We assume no slip at the solid-liquid interface and no stress at the liquid-air interface, so that:
\begin{eqnarray}
\label{CondLim}
v |_{z=0}&=&0\\
\partial_zv |_{z=h}&=&0\ .
\end{eqnarray}
The pressure being independent of the vertical coordinate $z$, Eq.~\eqref{Stokes} together with Eq.~\eqref{CondLim} lead to a Poiseuille flow of the form:
\begin{eqnarray}
v(x,z,t)&=&\frac{1}{2\eta}\,(z^2-2hz)\, \partial_x p\ .
\label{Poiseuille}
\end{eqnarray}
Conservation of volume can be expressed as:
\begin{eqnarray}
\partial_th&=&-\partial_x\int_0^h\textrm{d}z \,v\ .
\label{Mass}
\end{eqnarray}
Combining Eqs.~\eqref{Laplace}, \eqref{Poiseuille} and \eqref{Mass} leads to the two-dimensional capillary-driven thin film equation:
\begin{eqnarray}
\partial_th+\frac{\gamma}{3\eta}\,\partial_x\left( h^3\,\partial_x^{\,3}h \right)&=&0\ .
\label{TFE1}
\end{eqnarray}
Finally, Eq.~\eqref{TFE1} can be nondimensionalised by letting:
\begin{eqnarray}
h&=&H\, h_0\\
x&=&X\,h_0\\
t&=&T\,\frac{3\eta h_0}{\gamma}\ ,
\end{eqnarray}
where $h_0$ is the reference height at infinity. This leads to the dimensionless equation introduced above in Eq.~\eqref{TFE_AD}. This equation can be linearised by letting:
\begin{equation}
\label{defdelta}
H(X,T)=1+\Delta(X,T)\ , 
\end{equation}
and by assuming $\Delta\ll 1$:
\begin{eqnarray}
\partial_T\Delta+\partial_X^{\,\,4}\Delta&=&0\ .
\label{LTFE_AD}
\end{eqnarray}
The linear case corresponds to a situation in which the surface of a flat film is only slightly perturbed. Note that Eq.~\eqref{LTFE_AD} also describes surface diffusion phenomena leading to flattening of solid surfaces \cite{Mullins1958,Zhu2011}, kinetic growth \cite{Krug1993,Lai1991,Wolf1990}, or  grooving \cite{Mullins1957,Robertson1971}. Therefore, the following results apply to a broader class of physical situations.
\newline

%%%%%%%%%%%%%%%%%%%%%%%%%%%%%%%%%%%%%%%%%%%%%%%%%%%%%%%%%%%%%%%%%%%%%%

\section{Solving the linear equation}
In this part, we solve Eq.~\eqref{LTFE_AD} for $T>0$ and we characterise its solutions. The linear study is divided into four paragraphs. In the first one, we derive the Green's function and show that it is self-similar at all positive times. In the second one, we give the general formal solution and exploit it through two particular canonical examples. In the third one, we study the uniform convergence in time of the rescaled general solution towards the rescaled self-similar Green's function. In the fourth one, we conclude the discussion on the linear case with some general remarks.
\newline

%%%%%%%%%%%%%%%%%%%%%%%%%%%%%%%%%%%%%%%%%%%%%%%%%%%%%%%%%%%%%%%%%%%%%%

\subsection{Green's function and self-similarity}
Since Eq.~(\ref{LTFE_AD}) is a linear partial differential equation, it can be solved by calculating the Green's function $\mathcal G(X,T)$. This object is defined by:
\begin{eqnarray}
\left[\partial_T+\partial_X^{\,\,4}\right]\,\mathcal G(X,T)&=&\delta(X,T)\ ,
\label{Green1}
\end{eqnarray}
where $\delta$ denotes the Dirac distribution. Fourier transforms are defined as follows:
\begin{eqnarray}  \label{E:gp}
\mathcal G(X,T)&=&\frac{1}{(2\pi)^2}\int    \textrm{d}K    \textrm{d}\Omega\ \hat {\mathcal G}(K,\Omega)\,e^{i(\Omega T +KX)} \label{E:gp1}\\
\tilde{\mathcal G}(K,T)&=&\frac{1}{2\pi}\int    \textrm{d}\Omega\ \hat {\mathcal G}(K,\Omega)\,e^{i\Omega T}\label{E:gp2}\\
\hat {\mathcal G}(K,\Omega)&=&\int    \textrm{d}X    \textrm{d}T\ {\mathcal G}(X,T)\,e^{-i(\Omega T +KX)}\ .\label{E:gp3}
\end{eqnarray}
Taking the Fourier transform (see Eq.~(\ref{E:gp1}c)) of Eq.~(\ref{Green1}), and assuming that the Green's function vanishes at $X=\pm\infty$, consistent with the physical boundary conditions, leads to:
\begin{eqnarray}
\label{fourgreen}
\hat {\mathcal G}(K,\Omega)&=&\frac{1}{i\Omega +K^4}\ .
\end{eqnarray}
Using Eq.~(\ref{E:gp1}b) and Eq.~(\ref{fourgreen}), one obtains:
\begin{equation}
\tilde{\mathcal G}(K,T)=\textrm{Res}\left(\frac{e^{i\Omega T}}{\Omega -iK^4}\, ;\,iK^4  \right)\, \Theta(T)\ ,
\label{Caus}
\end{equation}
where Res\,$(f;z^*)$ denotes the complex residue of the function $f$ at $z=z^*$, and where $\Theta$ is the Heaviside function ensuring causality. Finally, expressing the residue and performing the inverse Fourier transform, defined in Eq.~(\ref{E:gp1}a), one obtains the Green's function: 
 \begin{eqnarray}
 \label{greenx}
\mathcal G(X,T)&=&\frac{1}{2\pi}\int    \textrm{d}K    \,     e^{-K^4T}            \,e^{iKX}\ ,
\end{eqnarray}
which is consistent with previous studies \cite{Bowen2006,Mullins1957}. Then, at finite time, let us change variables through:
\begin{eqnarray}X&=&UT^{1/4}\label{CDV1}\\ K&=&QT^{-1/4} \\ \mathcal G(X,T)&=&\breve{\mathcal{G}}(U,T)\ .
\label{CDV2}
\end{eqnarray}
This gives:
\begin{eqnarray}
\breve{\mathcal{G}}(U,T)
&=& \frac 1{T^{1/4}}\ \phi(U)\, \Theta(T)\ ,
\end{eqnarray}
where we introduced the auxiliary function:
\begin{equation}
\label{auxi}
\phi(U)=\frac{1}{2\pi}\,\int  \textrm{d}Q\,{e^{-Q^4}e^{iQU}}\ .
\end{equation}
\begin{figure}[t!]
\includegraphics[width=1\columnwidth]{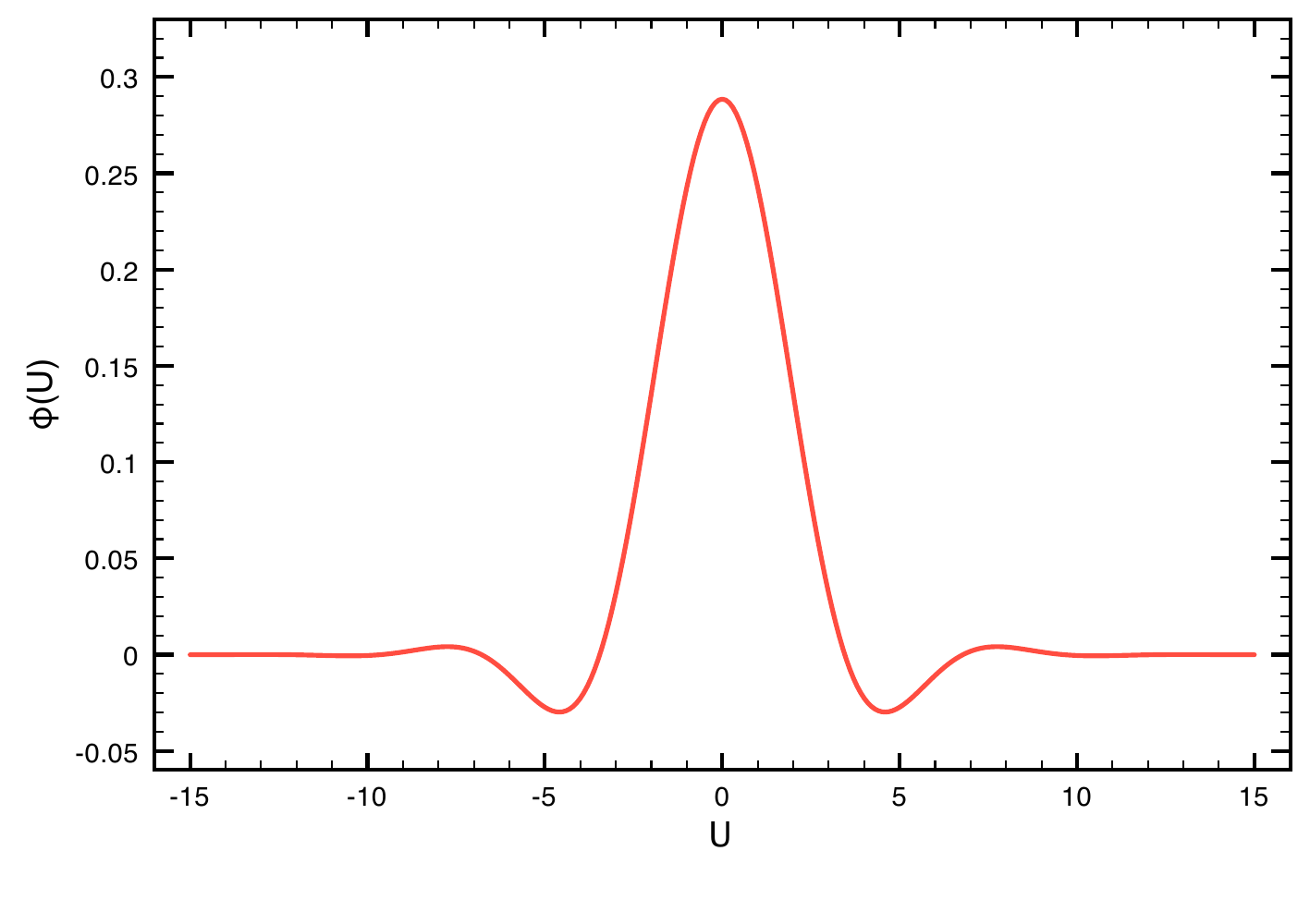}
\caption{Auxiliary function $\phi(U)={T^{1/4}}\ \breve{\mathcal G}(U,T)$, for positive times as given in Eq.~\eqref{HPG}, where $\breve{\mathcal G}(U,T)$ is the Green's function of the dimensionless linear two-dimensional capillary-driven thin film equation given in Eq.~\eqref{LTFE_AD}.}
\label{Fig2}
\end{figure}
The Green's function is thus \textit{self-similar of the first kind} \cite{Barenblatt1996}, at all positive times. Furthermore, the function $\phi$ is given for all $U\in \mathbb R$ by:
\begin{eqnarray}
\phi(U)&=& \frac1{\pi} \,\Gamma\left(\frac{5}{4}\right)\ _0H_{2}\left(\left\{ \frac12,\frac34\right\},\left(\frac{U}{4}\right)^4\right)\nonumber\\&-&\frac1{8\pi}\,U^2 \,\Gamma\left(\frac{3}{4}\right)\ _0H_{2}\left(\left\{ \frac54,\frac32\right\},\left(\frac{U}{4}\right)^4\right)
\label{HPG}
\end{eqnarray}
where the $(0,2)$-hypergeometric function is defined as \cite{Abramowitz1965,Gradshteyn1965}: 
\begin{eqnarray}
_0H_{2}\left(\left\{ a,b\right\},w\right)&=&\sum_{k\geq0} \frac{1}{(a)_k(b)_k}\,\frac{w^k}{k!}\ ,
\end{eqnarray}
with the Pochhammer notation $(.)_k$ for the rising factorial. The function $\phi$ is plotted in Fig.~\ref{Fig2}.
Other than the oscillatory behaviour which is directly related to the fourth spatial derivative of Eq.~(\ref{LTFE_AD}), this solution is qualitatively close to the point-source solution of the heat equation \cite{Fourier1822}, for which the same analytical treatment would lead to the well-known Green's function and to its specific self-similar variable $XT^{-1/2}$.
\newline

%%%%%%%%%%%%%%%%%%%%%%%%%%%%%%%%%%%%%%%%%%%%%%%%%%%%%%%%%%%%%%%%%%%%%%

\subsection{General solution}
\begin{figure}[t!]
\includegraphics[width=1\columnwidth]{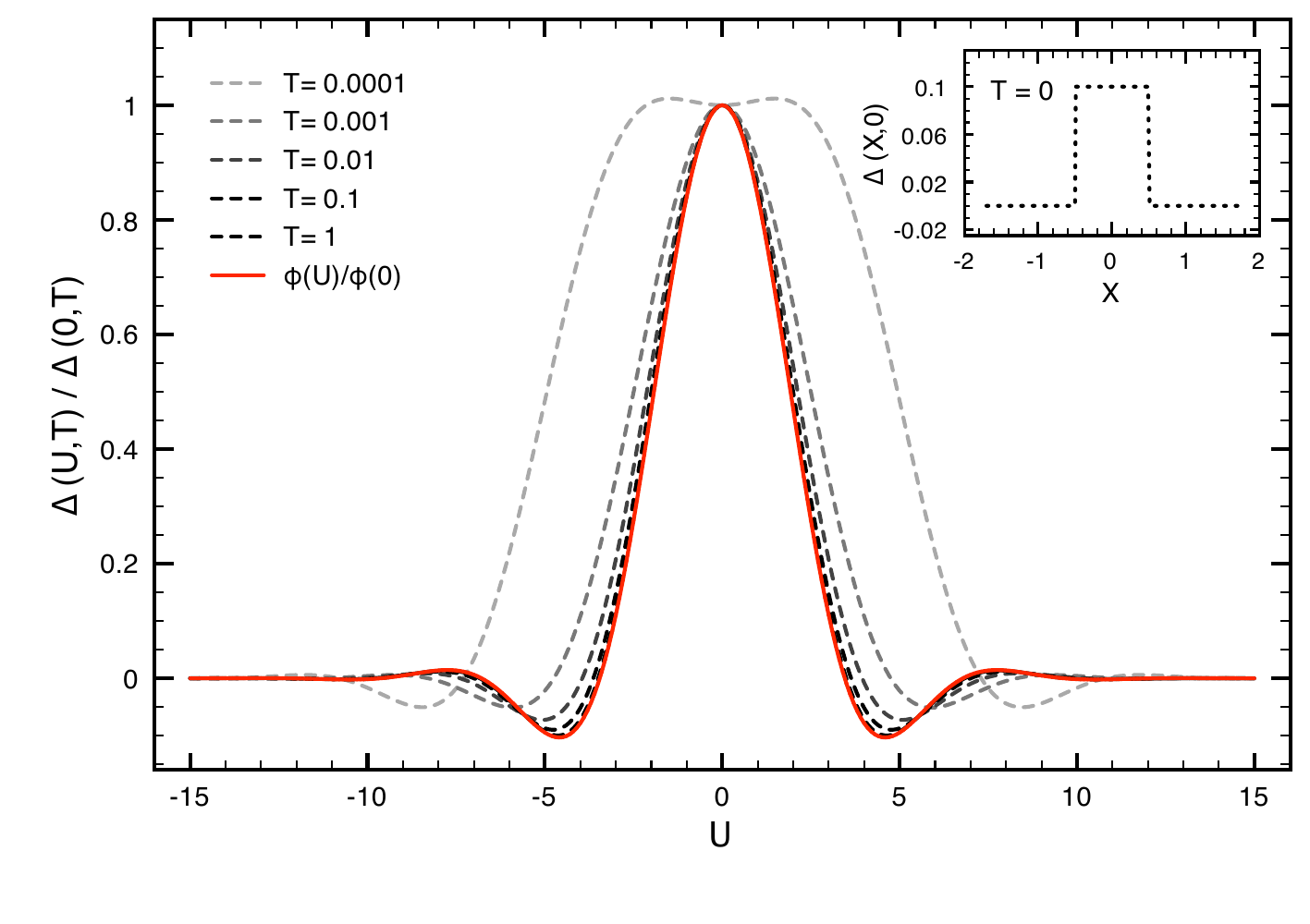}
\caption{Normalised analytical solution (see Eq.~\eqref{Gen_sol}) of the dimensionless linear two-dimensional capillary-driven thin film equation (see Eq.~\eqref{LTFE_AD}) as a function of the self-similar variable $U$ introduced in Eq.~\eqref{CDV1}. The initial profile is defined as a gate function of width unity, as shown in the inset. The solution is plotted for different dimensionless times and compared to its normalised asymptotic attractor given in Eq.~\eqref{HPG}.}
\label{Fig3}
\end{figure}

For a given summable \footnote{Along the present article, \textit{summable} means \textit{Lebesgue integrable}.\medskip} initial condition, $\Delta(X,0)=\Delta_0(X)$, the solution $\Delta(X,T)$ of Eq.~(\ref{LTFE_AD}) is given by the spatial convolution of $\Delta_0(X)$ with the Green's function:
  \begin{eqnarray}
\Delta(X,T)&=&\int \textrm{d}Y\,\mathcal G(X-Y,T)\,\Delta_0(Y)\ .
\label{Sol_xt}
\end{eqnarray}
Interestingly, Eq.~(\ref{Sol_xt}) implies that the Green's function of the problem is \textit{exactly} the point-source solution obtained from an initial Dirac spatial distribution: $\Delta_0(Y)=\delta(Y)$. Let us once again change variables through Eq.~(\ref{CDV1}). Then, Eq.~(\ref{greenx}) and Eq.~(\ref{Sol_xt}) lead to:
\begin{eqnarray}
\label{Gen_sol}
\Delta(X,T)&=&\breve{\Delta}(U,T)\\ =\frac{1}{2\pi\ T^{1/4}}&&\int \textrm{d}Q\,e^{-Q^4}\,e^{iQU}\int \textrm{d}Y\,e^{-iQY/T^{1/4}} \,\Delta_0(Y)\nonumber\ .
\end{eqnarray}
The particular case of an initial \textit{stepped film} \footnote{For comparison, in such a case one gets at $T>0$: 
$$\breve{\Delta}(U,T)=\frac{\theta_0}{2}\left(  1+\dashint  dQ \,\frac{1}{i\pi Q}\,e^{-Q^4}\,e^{iQU} \right)\ ,$$ 
where the dashed integral represents Cauchy's principal value and where $\theta_0$ is the amplitude of the step. This solution is self-similar for all $T>0$.} was studied in detail in a previous communication \cite{Salez2012}. In the present study, we have calculated the solutions for various summable initial conditions. We present two of them corresponding to canonical illustrations: a gate function (see Fig.~\ref{Fig3}) and a gaussian function (see Fig.~\ref{Fig4}). 
In each case, we plot the normalised analytical solution of Eq.~(\ref{LTFE_AD}) given in Eq.~\eqref{Gen_sol}, as a function of the self-similar variable $U$ introduced in Eq.~(\ref{CDV1}), the initial profile at $T=0$ being shown in the inset. The solution is plotted for different dimensionless times. For comparison, we also plot the normalised function $\phi(U)/\phi(0)$, as given in Eq.~\eqref{HPG}. As one can see, both profiles seem to converge in time towards the later  function that we shall henceforward call a \textit{universal self-similar attractor} (see Fig.~\ref{Fig1}). This convergence statement will be addressed in the next paragraph.
\newline

\subsection{Uniform convergence to the self-similar attractor}
For clarity, we shall restrict here to summable initial profiles that have non-zero algebraic volume, that is:
\begin{eqnarray}
\label{cond_mom}
\mathcal{M}_0\,\,=\int \textrm{d}X\,\Delta_0(X)&\,\,\neq\,\,&0\ .
\end{eqnarray}
The extension of the following results to the specific case of zero initial algebraic volume is understood as well and will be addressed in the next paragraph. For the time being, let us introduce the function:
\begin{eqnarray}
f(U,T)&=&\, \frac{T^{1/4}}{\mathcal{M}_0}\ \breve{\Delta}(U,T)\ .
\label{Pour_dev1}
\end{eqnarray}
\begin{figure}[t!]
\includegraphics[width=1\columnwidth]{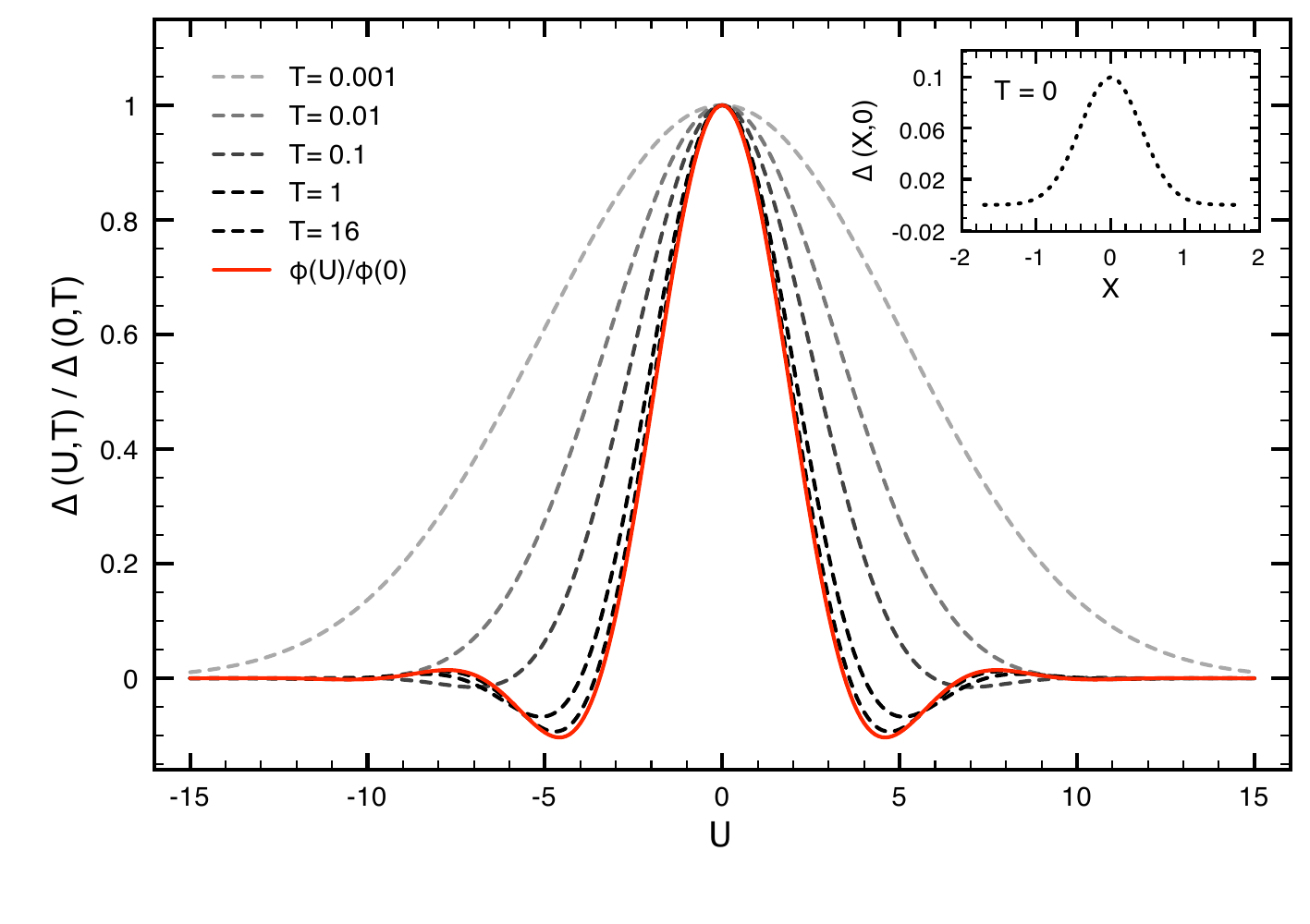}
\caption{Normalised analytical solution (see Eq.~\eqref{Gen_sol}) of the dimensionless linear two-dimensional capillary-driven thin film equation (see Eq.~\eqref{LTFE_AD}) as a function of the self-similar variable $U$ introduced in Eq.~\eqref{CDV1}. The initial profile is defined as a gaussian function of volume $0.1$, as shown in the inset. The solution is plotted for different dimensionless times and compared to its normalised asymptotic attractor given in Eq.~\eqref{HPG}.}
\label{Fig4}
\end{figure} 
According to Eq.~(\ref{auxi}), Eq.~(\ref{Gen_sol}) and Eq.~(\ref{Pour_dev1}), for all $T>0$ and for all $U\in \mathbb R$, one has:
\begin{equation}
\left| f(U,T)-\phi(U) \right| \leq \frac{a(T)}{2\pi|\mathcal{M}_0|}\ ,
\end{equation}
where we introduced:
\begin{equation}
a(T)= \int \textrm{d}Q\,e^{-Q^4}\int \textrm{d}Y\,\left|e^{-iQY/T^{1/4}}-1\right|\,\left|\Delta_0(Y) \right|\ .
\end{equation}
Therefore, one gets for all $T>0$:
\begin{eqnarray}
\left|\left|  f(U,T)-\phi(U) \,\right|\right|_{\infty,U} &\leq & \frac{a(T)}{2\pi|\mathcal{M}_0|}\ ,
\end{eqnarray}
where $|| ... ||_{\infty,U}$ is the uniform norm \footnote{$\|f\|_{\infty,x}=\sup\left\{\,\left|f(x,y)\right|,\,x\in \mathbb R\,\right\}. $\smallskip} with respect to the variable $U$. In order to conclude that the function $f$ is uniformly convergent \footnote{A function of two variables $f(x,y)$ is defined to be \textit{uniformly convergent} with respect to $y$ if $\displaystyle \lim_{y\rightarrow \infty}\|f\|_{\infty,x} = 0$.} in time towards $\phi$, it remains to show that: $\displaystyle \lim_{T\rightarrow \infty} a(T) =0$.
For this purpose, let us consider the auxiliary function defined by: 
\begin{eqnarray}
{m(Y,T)=\left(e^{-iQY/T^{1/4}}-1\right) \,\Delta_0(Y)}\ ,
\end{eqnarray}
for all $Y\in \mathbb R$, and for all $T>0$. This function naturally converges to the zero function when $T\rightarrow \infty$. In addition, for all $T>0$, and for all $Y\in \mathbb R$, one has:
\begin{eqnarray}
\label{condi1}
\left|m(Y,T)\right|&\leq& 2\,\left|\Delta_0(Y )\right|\ ,
\end{eqnarray}
where the right-hand side is a summable function. Invoking the continuity theorem for functions defined by a Lebesgue integral, one has:
\begin{eqnarray}
\label{lim1}
\lim_{T\rightarrow\infty} \int \textrm{d}Y\,m(Y,T)&=&0\ .
\end{eqnarray}
Similarly, we then consider the second auxiliary function defined by: 
\begin{eqnarray}
g(Q,T)=e^{-Q^4}\,\int \textrm{d}Y\,m(Y,T)\ ,
\end{eqnarray}
for all $Q\in \mathbb R$, and for all $T>0$. It converges to the zero function when $T\rightarrow \infty$, according to Eq.~(\ref{lim1}). In addition, using Eq.~(\ref{condi1}), for all $T>0$ and for all $Q\in \mathbb R$, one has:
\begin{eqnarray}
\left| g(Q,T) \right|&\leq& 2\,e^{-Q^4}\,\int \textrm{d}Y\,\left|\Delta_0(Y )\right|\ ,
\end{eqnarray}
where the right-hand side is a summable function. Once again, invoking the continuity theorem for functions defined by a Lebesgue integral, leads to:
\begin{eqnarray}
\lim_{T\rightarrow\infty}a(T)&=& 0\ . 
\end{eqnarray}
In summary, we demonstrated that the rescaled solution $f(U,T)$ of Eq.~(\ref{LTFE_AD}), given by Eq.~(\ref{Gen_sol}) and Eq.~(\ref{Pour_dev1}) for any summable initial profile, always converges \textit{uniformly} in time towards the universal self-similar attractor $\phi(U)$ defined in Eq.~(\ref{auxi}). In other words, we exhibited the intermediate asymptotics of the solutions of the linear capillary-driven thin film equation for flat boundary conditions and we demonstrated that it is simply given by the rescaled Green's function. We refer again to Fig.~\ref{Fig3} and Fig.~\ref{Fig4} for illustration of this result with a gate function and a gaussian function as initial profiles, respectively. The crucial point in these graphs is that the attractor is \textit{identical} for the two summable initial profiles, as summarised in Fig.~\ref{Fig1}. 
\newline

\subsection{Remarks} In order to conclude the discussion on the linear case, we enumerate five general remarks below.
\begin{itemize}
\item The latter calculations may be extended to all equations of the form:  
\begin{equation}
\left[\partial_T +\partial_X^{\,2m}\right]\Delta(X,T) =0\ , 
\end{equation}
with $m\in \mathbb N^*$. Indeed, other than the well known heat equation \cite{Fourier1822}, one can obtain the Green's function for higher even orders of the spatial derivative and extract analogous conclusions. However, the odd orders being free from dissipation are expected to lead to fundamentally different mathematical solutions.
\smallskip

\item According to our primary interest in thin films, we proved the uniform convergence of any summable solution towards the self-similar attractor by explicitly writing the Green's function of Eq.~(\ref{LTFE_AD}). However, this result is more general. In fact, let us consider any diffusive-like linear partial differential equation of two variables for which the Green's function $G(X,T)$ is self-similar of the form:
\begin{equation}
G(X,T)=\frac{1}{T^\beta}\ \phi\left(\frac{X}{T^{\alpha}}\right)\ ,
\end{equation}
where $\alpha$ and $\beta$ are strictly positive real numbers, and where $\phi$ is a bounded function on $\mathbb R$. Then, the uniform convergence is straightforward to demonstrate, as soon as Eq.~(\ref{Sol_xt}) and Eq.~(\ref{cond_mom}) remain satisfied, with a summable initial profile.
\smallskip

\item Another interesting feature is the following. Let us call $\Delta_{[\Delta_0]}(X,T)$ the solution of Eq.~(\ref{LTFE_AD}), for an initial profile $\Delta_0(X)$. Then, the evolution $ \Delta_{[\Delta_0']}(X,T)$ of the first derivative $\Delta_0'$ of the previous initial condition is simply given by $\partial_X  \Delta_{[\Delta_0]}(X,T)$. In terms of the self-similar variable introduced in Eq.~(\ref{CDV1}), one gets:
\begin{equation}
\breve{\Delta}_{[\Delta_0']}(U,T)=\frac{1}{T^{1/4}}\,\partial_U  \breve{\Delta}_{[\Delta_0]}(U,T)\ .
\end{equation}
This notably makes the link with our previous work on the particular case of an initial stepped film \cite{Salez2012}. Indeed, the Dirac distribution being the first derivative of the Heaviside's distribution, the Green's function is thus simply given by the derivative of the solution of Eq.~(\ref{LTFE_AD}) obtained for a stepped initial condition. Generalised to higher order derivatives, this relation naturally becomes:
\begin{equation}
\breve{\Delta}_{[\Delta_0^{(n)}]}(U,T)=\left[\frac{1}{T^{1/4}}\,\partial_U \right]^{n} \breve{\Delta}_{[\Delta_0]}(U,T)\ .
\end{equation}
\smallskip

\item One may as well wonder what happens in the particular case of zero algebraic volume \footnote{One could for instance imagine a \textit{dip} followed by a \textit{bump} of identical shape, as a pathological initial profile.}, that is when $\mathcal{M}_0=0$. The answer is that there is still an attractive self-similar regime given by the first non-zero derivative $\phi^{(n)}(U)$, under the condition that this quantity is summable.
\smallskip

\item At last, we have seen that for any summable initial condition there is long-term self-similarity of the general solution of Eq.~(\ref{LTFE_AD}). The question arrises to know whether there exists some specific summable initial conditions that generate solutions that are self-similar \textit{at all positive times}. Interestingly, the solution obtained from the Heaviside initial condition was shown to be self-similar at all positive times for other boundary limits \cite{Salez2012}. In the present case, let us impose the following constraint at all times $T>0$:
\begin{equation}
\label{defsim}
\breve{\Delta}(U,T)=T^{\alpha/4}F(U)\ ,
\end{equation}
with $U$ as defined in Eq.~(\ref{CDV1}). By changing variables in Eq.~(\ref{Gen_sol}), it is straightforward to see that the initial profile must be homogeneous of degree $\alpha$, meaning that for all real numbers $k$ and $Y$, one has:
\begin{equation}
\label{hom}
\Delta_0(kY)=k^{\alpha}\Delta_0(Y)\ .
\end{equation}
Finally, using Eq.~(\ref{cond_mom}) and Eq.~(\ref{hom}), the algebraic volume $V$ satisfies:
\begin{eqnarray}
V&=&\int \textrm{d}X\ \Delta(X,T)\\
&=&T^{(\alpha+1)/4}\int \textrm{d}U\,F(U)\ ,
\end{eqnarray}
which implies that $\alpha=-1$, since $V=\mathcal{M}_0$ by volume conservation. Therefore, the only initial profile with finite algrebraic volume that exhibits the self-similarity of Eq.~(\ref{defsim}) at all positive times is the Dirac distribution. Recalling the remark made after Eq.~(\ref{Sol_xt}), this means that the Green's function is the only summable solution that is self-similar at all positive times, with the definition of Eq.~(\ref{defsim}). 
\newline
\end{itemize}

%%%%%%%%%%%%%%%%%%%%%%%%%%%%%%%%%%%%%%%%%%%%%%%%%%%%%%%%%%%%%%%%%%%%%%

\section{Extension to the nonlinear equation}
\begin{figure}[t!]
\includegraphics[width=1\columnwidth]{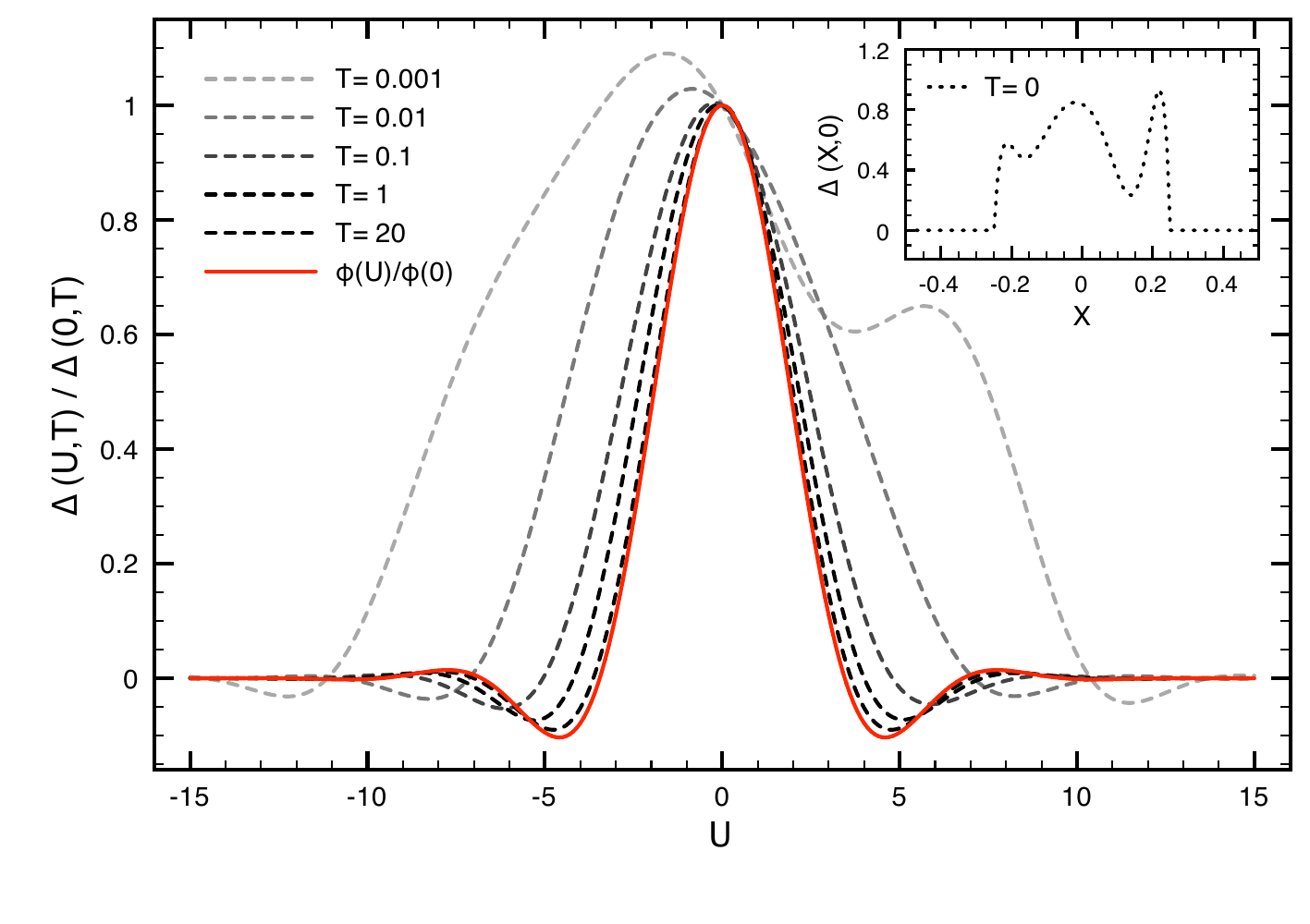}
\caption{Normalised numerical solution of the dimensionless nonlinear two-dimensional thin film equation given in Eq.~\eqref{TFE_AD}, as a function of the self-similar variable $U$ defined in Eq.~\eqref{CDV1}, for an initial profile given by an arbitrary function with compact support, as shown in the inset. The solution is plotted for different dimensionless times and compared to the normalised universal asymptotic attractor of the linear case given in Eq.~\eqref{HPG}.}
\label{Fig5}
\end{figure} 
In this last part, we extend the previous results to the nonlinear case through a numerical scheme. The excess profile $\Delta(X,T)$ is still defined by Eq.~(\ref{defdelta}), but without any restriction on its amplitude. Thus, we consider the full nonlinear partial differential equation given in Eq.~\eqref{TFE_AD}. This equation has not been solved analytically yet, but we recently solved it numerically in various geometries \cite{Salez2012a}. The numerical procedure we used is a finite difference method developed in \cite{Bertozzi1998,Zhornitskaya2000}. It ensures capillary energy and entropy dissipation, as required from \cite{Bernis1990}. In addition to volume conservation, it has been shown that this method ensures positivity of the height profile $H(X,T)$ \cite{Zhornitskaya2000}. Using this numerical scheme, we verified the existence of an attractive self-similar regime for several \textit{arbitrary} initial profiles, with \textit{compact support} as required from the algorithm. For instance, we plot in Fig.~\ref{Fig5} the normalised numerical solution of Eq.~\eqref{TFE_AD} as a function of the self-similar variable $U$ defined in Eq.~(\ref{CDV1}), the initial profile being shown in the inset. The solution is plotted for different dimensionless times. We also plot for comparison the normalised function $\phi(U)/\phi(0)$, as given in Eq.~\eqref{HPG}. As we see, the latter is thus an attractor for any initial profile with compact support (see Fig.~\ref{Fig1}). What is remarkable here is that the attractor of the \textit{nonlinear} case is actually equal to the rescaled Green's function obtained analytically in the \textit{linear} case. Even though this may seem unexpected, it turns out to be quite natural when one realises that as time goes, the initial profile progressively collapses towards the flat film equilibrium shape, thus bringing Eq.~\eqref{TFE_AD} closer to its linearised form of Eq.~\eqref{LTFE_AD}.
\newline

%%%%%%%%%%%%%%%%%%%%%%%%%%%%%%%%%%%%%%%%%%%%%%%%%%%%%%%%%%%%%%%%%%%%%%

\section*{Conclusion}
We reported on the intermediate asymptotics of the two-dimensional capillary-driven thin film equation for an arbitrary summable initial profile. First, we derived an analytical solution of the linearised equation. The solution was obtained by seeking the Green's function, which was found to be given by a combination of generalised hypergeometric functions. As schematised in Fig.~\ref{Fig1}, we then proved that any summable initial condition leads to the uniform convergence in time of the rescaled solution towards a universal self-similar attractor that is proportional to the rescaled Green's function of the problem; the proportionality factor being equal to the initial algebraic volume. This result appears to be a more general result in diffusive-like processes that are characterised by a self-similar Green's function with decaying amplitude. At last, we were able to conjecture from compact-support numerical results, as well as to justify on a physical basis, the extension of this convergence behaviour to the nonlinear equation. The important outcome is that the universal self-similar attractor of the nonlinear case appears to be precisely the one of the linear case, that is the rescaled Green's function. 
\smallskip

The recent excellent agreements between thin film theories and experiments with stepped polymer films \cite{McGraw2011,McGraw2012,Salez2012a,Salez2012}, and polymer droplets on identical films \cite{Cormier2012,Salez2012a}, are very encouraging for the physical relevance of the present analysis. Thus, experimental implementation with viscous nanofilms should be performed in near future. This theoretical work may also be extended to other thin film equations \cite{Bertozzi1998}, such as the one describing flows in Hele-Shaw cells \cite{Barenblatt1996} and the one governing gravity-driven flows \cite{Decre2003,Huppert1982,Kondic2003}. One could also address coalescence phenomena \cite{Hernandez2012,Ristenpart2006} and various diffusive processes, in an identical way. Moreover, it should be feasible to exhibit the necessary conditions acting on a given partial differential equation for showing such a uniform convergence behaviour. Finally, characterising further the conditions for finite-time convergence may be of fundamental importance for the field of intermediate asymptotics, and thus for nonlinear physics in general. 
\newline

%%%%%%%%%%%%%%%%%%%%%%%%%%%%%%%%%%%%%%%%%%%%%%%%%%%%%%%%%%%%%%%%%%%%%%

\section*{Acknowledgments} 
The authors warmly thank  Kari Dalnoki-Veress, Joshua D. McGraw and Oliver B\"aumchen for a very stimulating ongoing collaboration. They are also grateful to Ken Sekimoto, Mark Ediger and Alexandre Darmon for interesting discussions. Finally, they thank the \'Ecole Normale Sup\'erieure of Paris and the Fondation Langlois for financial support.

%%%%%%%%%%%%%%%%%%%%%%%%%%%%%%%%%%%%%%%%%%%%%%%%%%%%%%%%%%%%%%%%%%%%%%

\end{document}